
\def\date{\ifcase\month\or January\or February \or March\or April\or May %
\or June\or July\or August\or September\or October\or November           %
\or December\fi\space\number\day, \number\year}                          %
\magnification=\magstep1   
\hsize=15truecm            
\vsize=20truecm            
\parindent=0pt             
\line{\hfill\date}         
\vskip2\baselineskip       
\catcode`\@=11             
\newfam\msbfam             
\newfam\euffam             
\font\tenmsb=msbm10        
\font\sevenmsb=msbm7       
\font\teneuf=eufm10        
\font\seveneuf=eufm7       
\textfont\msbfam=\tenmsb   
\scriptfont\msbfam=\sevenmsb%
\textfont\euffam=\teneuf   
\scriptfont\euffam=
\seveneuf                  
\def\N{{\fam\msbfam N}}    
\def\Q{{\fam\msbfam Q}}    
\def\Z{{\fam\msbfam Z}}    
\def\cc{{\fam\euffam c}}    
\def\U{{\fam\euffam U}}    
\catcode`\@=12             
\def\.{{\cdot}}            
\def\<{\langle}            
\def\>{\rangle}            
\def\({\big(}              
\def\){\big)}              
\def\hat{\widehat}         
\def\card{\mathop{\rm card}
\nolimits}                 
\def\defi{\buildrel\rm def 
\over=}                    
\def\ssk{\smallskip}       
\def\msk{\medskip}         
\def\bsk{\bigskip}         
\def\nin{\noindent}        
\def\cen{\centerline}      
\def\implies
{\hbox{$\Rightarrow$}}     
\def\iff
{\hbox{$\Leftrightarrow$}} 
\font\smc=cmcsc10          
\def\subsetneq{\;\rlap{\raise3pt\hbox{${\subset}$}}%
{\raise-3pt\hbox{${\neq}$}}\;} 
\def\2{{\bf2}}             
\def\0{{\bf0}}             
\def\1{{\bf1}}             
\def\lead{\leaders\hbox to 1.5ex{\hss${.}$\hss}\hfill}
\def\arr{\hbox to 30pt{\rightarrowfill}}
\def\larr{\hbox to 30pt{\leftarrowfill}}

\long\def\alert#1{\parindent2em\smallskip\line{\hskip\parindent\vrule%
\vbox{\advance\hsize-2\parindent\hrule\smallskip\parindent.4\parindent%
\narrower\noindent#1\smallskip\hrule}\vrule\hfill}\smallskip\parindent0pt}

\def\bx#1{\lower4pt\hbox{
\vbox{\offinterlineskip
\hrule
\hbox{\vrule\strut\hskip1ex\hfil{\uvertspace #1}\hfill\hskip1ex}
\hrule}\vrule}}
\newcount\litter                                         
\def\newlitter#1.{\advance\litter by 1                   
\edef#1{\number\litter}}                                 
\def\qedbox{\hbox{$\rlap{$\sqcap$}\sqcup$}}              
\def\qed{\nobreak\hfill\penalty250 
\nobreak\hfill\qedbox\vskip1\baselineskip\rm}            

\newlitter\bourb.
\newlitter\bour.
\newlitter\comf.
\newlitter\ferher.
\newlitter\gart.
\newlitter\gelb.
\newlitter\cont.
\newlitter\hart.
\newlitter\hern.
\newlitter\HHM.
\newlitter\hewross.
\newlitter\hof.
\newlitter\layer.
\newlitter\compbook.
\newlitter\kall.
\newlitter\khara.
\newlitter\klepp.
\newlitter\mont.
\newlitter\segal.
\newlitter\pete.
\newlitter\ribes. 
\newlitter\ross.
\newlitter\rudina.
\newlitter\rudinb.
\newlitter\saeki.
\newlitter\saxl.
\newlitter\smith.
\newlitter\thomas. 
\newlitter\zele.
\newlitter\zel.


\cen{\bf Nonmeasurable subgroups of compact groups}
\msk
\cen{by Salvador Hern\'andez, Karl H. Hofmann, and Sidney A. Morris}
\bsk

{\bf Abstract.}\quad 

In 1985 S.~Saeki and K.~Stromberg published the following  question:
{\it Does every infinite compact group have
a subgroup which is not Haar measurable?}
An affirmative answer is given for all compact groups with
the exception of some metric profinite groups known as
strongly complete. In this
spirit it is also shown that every compact group contains a 
non-Borel subgroup. 

\bigskip 

\cen{\bf Introduction}

\msk

For a compact group $G$,
{\it measurable} means ``measurable with respect to
the unique normalized Haar measure $\mu$ on $G$.'' Since Haar measure is a
Borel measure, every Borel subset of $G$ is measurable. A subset
$S\subseteq G$ is a {\it  null set} if $\mu(K)=0$ for
each compact subset $K$ of $S$, and if for each $\epsilon>0$ there is an
open neighborhood of $S$ such that $\mu(U)<\epsilon$.
Every subset of a  null set  is
measurable (see [\bour], Paragraph after Chap IV, \S 5, n$^\circ$ 2,
Definition 3, p. 172, or [\hewross], p.~125, Theorem 11.30).

\msk
The topic of subsets of a (locally) compact group which are not measurable
with respect to Haar measure 
is a wide field. Hewitt and Ross provide an instructive and far-reaching
discussion of this topic in [\hewross], pp.\ 226ff.
The present question differs insofar as in this paper
we are looking for the existence of  nonmeasurable
subgroups rather than just nonmeasurable subsets.
\msk
{\bf Question 1.} [\saeki] \quad {\it Does every infinite compact group contain a 
nonmeasurable subgroup?}
\bsk
For abelian compact groups 
{\smc Comfort} et al [\comf] showed the existence of nonmeasurable subgroups. 
See also {\smc Kharazishvili} [\khara]. For some
partial answers in the noncommutative case see Gelbaum [\gelb], 4.45.
\bsk
\eject

\cen{\bf 1. The standard background material}
\msk
We now present a systematic approach towards answering Question 1.
First, we introduce some pertinent notation. 

The terms $|G|$ 
and $\card(G)$ equivalently denote  the cardinality of $G$, notably
when $G$ is a group. Secondly if $H$ is a subgroup of $G$ and $G/H$
is the set of cosets $gH$, $g\in G$, then $(G:H)\defi \card(G/H)$
denotes the {\it index of $H$ in $G$}. 
We say that
{\it $H$ has countably infinite index (in $G$)} iff $(G:H)=\aleph_0$,
and that
{\it $H$ has countable index (in $G$)} iff $\card(G/H)\le\aleph_0$.

\msk\bf
Proposition 1.1. \quad \it Let $G$ be a compact group and  $H$
a subgroup.

{\rm(a)} If $H$ has countably infinite index, 
then $H$ is nonmeasurable. In particular, $H$ is not a Borel subset.

\ssk

{\rm(b)} If $H$ is measurable, then either it has measure $0$, or it has
positive measure in which case it is open (thus having finite
index).
\ssk

\ssk

{\rm(c)} If $H$ has finite index in $G$ and is not closed, then
$H$ is nonmeasurable. In particular, a countable index subgroup $H$
of $G$ is either closed with finite index or is nonmeasurable.

\ssk

{\rm(d)} If $H$ is nonmeasurable in $G$,
then $\overline H$ is an open
(and therefore finite index) subgroup of $G$.

\ssk

{\rm(e)} If $H$ is a finite index subgroup,
$G=H\cup g_1H\cup\cdots\cup g_nH$, then the largest
normal subgroup  $N=H \cap g_1Hg_1^{-1}\cap\cdots\cap g_nHg_n^{-1}$
has finite index in $G$.

\ssk

{\rm(f)} Assume that $H$ is nonmeasurable and that $N$ is the largest
normal subgroup contained in the open subgroup $\overline H$.
Then $N$ is open and $H\cap N$ is dense in $N$ and nonmeasurable
in $N$.

{\rm(g)} Assume that $f\colon G\to G_1$ is a surjective morphism of
compact groups and that $H_1\subseteq G_1$ is a nonmeasurable 
subgroup of countable index. 
Then $H\defi f^{-1}(H_1)$ is a nonmeasurable subgroup of $G$.

\msk\bf
Proof.\quad \rm (a) (See [\hewross], p.\ 227 and [\saeki], Remark on p.\ 373.) 
Let $\{g_1=1,g_2,\dots\}$
be a system of representatives for $G/H$, that is
$$G=\bigcup_{n=1}^\infty g_nH,\quad \hbox{a disjoint union.}\leqno(1)$$
Suppose that $H$ is measurable. Then $g_nH$ is measurable for all $n$
and $\mu(g_nH)=\mu(H)$ by the invariance of Haar measure.
So (1) implies
$$1=\mu(G)=\sum_{n=1}^\infty\mu(g_nH)=
\sup_{N=1,2,\dots}\sum_{n=1}^N\mu(g_nH)= \sup_{N=1,2,\dots}N\.\mu(H).\leqno(2)$$
In particular, $\{N\.\mu(H):N=1,2,\dots\}$ is a bounded set of nonnegative
numbers, and this implies $\mu(H)=0$. Then $N\.\mu(H)=0$ for all
$N=1,2,\dots$ and so $\sup_{N=1,2,\dots}N\.\mu(H)=0$. This contradicts (2)
and therefore our supposition must be false. That is, $H$ is nonmeasurable.

\msk

(b) If $H$ is measurable and has positive measure,
then by [\hewross], 20.17, p.\ 296, the group
$H=HH$ has inner points, and thus is open.

\msk

(c) If $H$ has finite index in $G$, as in (a) above, let
  $\{g_1=1,g_2,\dots,g_N\}$ be a system of representatives
for $G/H$. Assume that $H$ is measurable. Then $1=\mu(G)=
\sum_{n=1}^N \mu(g_nH)=N\.\mu(H)$. Thus $\mu(H)={1\over N}>0$.
Then $H$ is an open subgroup by (b) and thus is also closed.

\msk

(d) Indeed, $H$ is nowhere dense iff $\overline H$ has no inner points.
So $\mu(H)=0$. Then $H$ is a subset of a null set and therefore
is measurable (and is a null set).

\msk

(e) This is straightforward.

\msk

(f) $N$ is open and of finite index by (d) and (e).
Now $N=\overline{N\cap H}$
since $N$ is open and $H$ is dense in $\overline H$. Also, since $N$
is an identity neighborhood, $\overline H\subseteq HN \subseteq \overline H$,
that is, $\overline H= HN=NH$. There are elements $h_1,\dots,h_k\in H$
such that $\overline H=N\cup h_1N\cup\cdots\cup h_kN$ is a coset decomposition
of $\overline H$.

\msk

Now suppose that $N\cap H$ is measurable.
Then $h_k(N\cap H)=h_kN\cap H$ is measurable
for all $k=1,\dots,m$, and so 
$(N\cap H)\cup (h_1N\cap H)\cup\cdots\cup(h_kN\cap H)=
\overline H \cap H=H$ 
is measurable in contradiction to the hypothesis on $H$. 
If $\mu$ is Haar measure on $G$, then
$(G:N)^{-1}\mu|N$ is normalized Haar measure on $N$. 

\msk

(g) If the index of $H_1$ in $G_1$ is infinite, then
$|G/H|$ is infinite, whence  $H$ is nonmeasurable by (a) above.
If  $(G_1:H_1)<\infty$, then $(G:H)<\infty$ since
 $G/H\cong G_1/H_1$. If $H$ were measurable,
then it would be open in $G$ by (b) and thus $H_1$ would be open
in $G_1$ which is not the case. So $H$ is nonmeasurable. \qed

Regarding condition (c) above we should note right away
that an infinite algebraically simple compact group such as
SO(3) (see [\compbook], Theorem 9.90)
does not contain any proper finite index subgroups 
in view of (e) while,
as we shall argue in Theorem 2.3 below, it does contain
countably infinite index subgroups.
On the other hand, a power $A_5^\N$ with the alternating
group $A_5$ of 60 elements does not contain
any countably infinite index subgroup as {\smc Thomas} shows in
Theorem 1.10 of [\thomas], while it does contain
nonclosed proper finite index subgroups.
In [\klepp], {\smc Kleppner} shows (in terms of homomorphisms onto finite
discrete groups) that nonopen finite index normal subgroups are
nonmeasurable. 

\msk
In order to better understand the focus of our observations let us say
that we may distinguish the following classes of compact groups:

Class 1:  compact groups having subgroups of countably infinite index;

Class 2:  compact groups having nonclosed subgroups of finite index;

Class 3:  compact groups having in which every countable index subgroup
is open closed.

\msk
In the direction of answering Question 1, the listing of known facts in 
Proposition 1.1 allows us to say that all groups in Classes 1 and 2 have
nonmeasurable subgroups. The group SO$(3)$ is a member of Class 1 but is not
in Class 2. The group  $A_5^\N$ belongs to Class 2 and not to Class 1.
If $G_1$ is a member of Class 1 and $G_2$ is a member of Class 2, then
$G_1\times G_2$ is a member of the intersection of the two classes.
   In the end we have to focus on Class 3, the complement of the union of
the first two classes; however, it will serve a useful purpose to 
understand how big this union is and where familiar categories of compact
groups are classified in this system.

\msk
The following discussion of  examples show  how members
of Class 1 and 2 may arise.
 For this purpose let $K$ be an arbitrary compact nonsingleton group.
Let $X$ be an infinite set endowed with the
discrete topology, for instance
$X=\N$. Then the compact group $G=K^X$ has an alternative description.
Indeed we consider $X$ as a subset of its Stone \v Cech compactification
$\beta X$ and note that, due to the compactness of $K$, every element
$f\in G$, that is, every function $f\colon X\to K$ has a unique extension
to a continuous function $\overline f\colon \beta X\to K$. The 
function $f\mapsto\overline f:G\to C(\beta X,K)$ is an isomorphism of
groups if we give $C(\beta X, K)$ the pointwise group operations.
If we endow $C(\beta X,K)$ with the topology of pointwise convergence on
the points of $X$, then $C(\beta X,K)$ is a compact group and 
$f\mapsto \overline f$ is an isomorphism of compact groups with the inverse
$F\mapsto F|X$. We shall identify $G$ and $C(\beta X, K)$ and note that 
$G$ has a much finer topology, namely, that of uniform convergence
on compact subsets of $\beta Y$ (to compact-open topology)
giving us a topological group $\Gamma$ with the same
underlying group as $G$. Now let $y\in \beta X$
and let $H$ denote a proper closed subgroup of $K$.
Then $G_{y,H}=\{f\in \Gamma: f(y)\in H\}$ is a closed subgroup of $\Gamma$.
 We record the following lemma:
\msk

\bf Lemma 1.2. \it The following statements are equivalent:
{\parindent1.5em 

\item{\rm(1)} $G_{y,H}$ is closed in $G$.

\item{\rm(2)} $y\in X$.

}
\msk

\bf Proof.\quad \rm   Since $y\in X$ implies that, for the
continuous projection  $p_y:G\to K$, $p_y(f)=f(y)$, the set $G_{y,H}$ is just
$p_y^{-1}(H)$, (1) follows trivially from (2).

Now suppose (1) is true and (2) is false. 
We must derive a contradiction.  
Since $\beta X$ is zero dimensional (in fact extremally
disconnected), the point $y$ has a basis $\U$ of open-closed neighborhoods
Let $g\in K\setminus H$ and define, for each $U\in\U$, a continuous function
 $f_U\colon \beta X\to K$ by 
$$f_U(z)=\cases{g &if $z\in\beta X\setminus U$,\cr
               1 &if $z\in U$.\cr}\leqno(*)$$
Since $G$ is compact, there is a cofinal net $(U_j)_{j\in J}$ in $\U$
such that $f=\lim_{j\in J} f_{U_j}$ exists in $G$. By $(*)$ we have
$f\in G_{h,H}$. Since $y$ is not isolated in $\beta X$ as (2) fails,
 also from $(*)$ we have a net $x_j\in X\setminus U_j$
converging to $y$ such that $f_{U_j}(x_j)=g$. Let $N$ be an open neighborhood
of $g$ in $K$ with $1\notin N$. 
 The continuity of $f$ implies the existence of a closed neighborhood
 $W$ of $y$ in $\beta X$ such that  $f(W)\cap N=\emptyset$.
Now let $k\in J$ be such that $j\ge k$ implies $x_j\in W$. 
Then for $i\ge k$ we have 
$g =\lim_{j\ge k}f_{U_j}(x_i)=f(x_i)\in f(W)$, and thus
$g\notin N$, a contradiction.\qed  
\ssk
The argument shows in fact that 
for $y\in\beta X\setminus X$, the proper
subgroup $G_{y,H}$ is dense in $G$.
\ssk

\msk
\bf
\nin
Corollary 1.3. \quad\it Every compact group $G$ of the form $G=K^X$, for an
infinite set $X$ and a profinite group $K$, has nonmeasurable subgroups.
\msk
\nin
\bf Proof.\quad \rm Let $H$ be a proper subgroup of $K$ of finite
index. Then for each $y\in \beta X\setminus X$, the subgroup
$G_{y,H}$ fails to be  closed by Lemma 1.2. 
On the other hand, since $G/G_{y,H}\cong K/H$ algebraically,
$G_{y,H}$ has finite index   and thus is
not measurable by 1.1(c).\qed

In particular, 
\ssk{\it
for each finite group $F$ and each infinite set $X$,
the profinite group $F^X$ has nonmeasurable subgroups.}

\bsk
As a an extension of Corollary 1.3 we mention the following observation
which is obtained as a simple application of Proposition 1.1(g).
\msk

\bf
Corollary 1.4.\quad \it If $G=\prod_{j\in J} G_j$ for a family of
compact groups $G_j$ and there is an infinite subset $I\subseteq J$
such that $G_j\cong K$ for all $j\in I$ and for a profinite group $K$,
then $G$ has nonclosed finite index and therefore  
nonmeasurable subgroups. \bf

\msk

Proof. \quad \rm We may identify each $G_j$ with $K$ for $j\in I$
and define the morphism  $f\colon G\to K^I$ as the obvious partial 
product. Then $K^I$ has a nonmeasurable subgroup $H_1$ 
by Corollary 1.3. So $H\defi f^{-1}(H_1)$ is a nonmeasurable subgroup of $G$
by Proposition 1.1(g).\qed

\eject

\cen{\bf 2.  The case of countably infinite index subgroups}

\msk

In [\HHM], Corollary 1.2 we noted that
every uncountable  abelian group has a proper subgroup $H$ of
countable index. (See also [\hewross], p.~227.)
Accordingly, by Proposition 1.1(a),  we have

\msk
\bf
Proposition 2.1. [\comf]
 {\it An infinite compact abelian group has
a nonmeasurable subgroup}. \qed
 
\msk
\rm

A bit more generally, we have
the following observation as an immediate consequence of Proposition 2.1
and Proposition 1.1(c).

\msk
\bf
Corollary 2.2.\quad \it If  the algebraic commutator group $G'$ of
a compact group $G$ has infinite
index in $G$, then $G$ has nonmeasurable subgroups. 
If the subgroup $G'$ has finite index it is either 
open closed or nonmeasurable. \qed

\msk
Here is a partial answer to Question 1:

\msk
\bf
Theorem  2.3. \quad\it Every infinite compact  group $G$ that is not profinite
has a subgroup of countable index and thus contains a nonmeasurable
subgroup.

\msk
\bf
Proof.
\quad\rm
Assume that $G$ is not profinite. Let $N$ be a closed normal subgroup
of $G$ such that $C=G/N$ is an infinite compact Lie group ([\compbook],
Corollary 2.43).

We will show that  $C$  contains
a subgroup with countable index; for
the  pullback to $G$ of a subgroup of countable index
in $C$ yields a subgroup of countable index in $G$.
A subgroup of the identity component $C_0$ of $C$ 
with countable index has countable index
in $C$ as $C/C_0$ is finite. Now $C_0$ is a compact connected
Lie group and we claim that it has  a subgroup of countable index.

The commutator subgroup $C_0'$ of $C_0$ is closed 
([\compbook], Theorem 6.11)
and so, if $C_0'\ne C_0$, then $C_0/C_0'$ is a connected
abelian Lie group and thus is infinite and therefore contains
a subgroup of countable index by Lemma 1.2.
Thus $C_0$ has a countable index subgroup.

Next we assume that $C_0'=C_0$ and thus that $C_0$
is semisimple and there is a homomorphism onto a
compact connected simple and centerfree Lie group  $S$.
(See [\compbook], Theorem 6.18.)
Now $S$ has no subgroup of finite index, because  if $H$
were a finite index subgroup of $S$, then the intersection of
the finitely many conjugates of $H$ would be a finite index normal
subgroup which cannot exist (see [\compbook], Theorem 9.90).

Now $S$ has a faithful linear representation
as an orthogonal matrix group (cf.\ [\compbook]. Corollary 2.40).
By a theorem of {\smc Kallman} [\kall], therefore $S$ has a faithful
algebraic representation as a permutation group on $\N$.
Since all orbits of $S$ on $\N$ are countable or finite,
the isotropy groups  all have countable or finite index. Since
$S$ has no finite index subgroups, it does have a countable
index subgroup which finally pulls back to a finite index
subgroup of $G$. \qed
\bsk

At this stage there remains the case of  profinite groups.
\ssk

Before we address this case let us observe that even 
in the compact abelian case
the issue of countably infinite index subgroups is far from trivial.
{}From [\hof] we quote  (cf.\ also [\compbook], Theorem 8.99):

\msk
\bf
Theorem 2.4. \it\quad There is a model of set theory
in which there is a compact
group $G$ with weight $\aleph_1=2^{\aleph_0}$ such
that the arc component factor group
$\pi_0(G)=G/G_a$ is algebraically isomorphic to $\Q$. In particular, the
arc component of  $G_a$ of the identity is a countably
infinite index subgroup and, accordingly, is nonmeasurable.

\rm\bsk

\cen{\bf 3. Profinite groups}
\msk

We record next that not all compact groups have countably infinite index
subgroups:
\msk
{\bf Example 3.1.}\rm \quad Let $A_5$ the alternating group on 
five elements, the smallest finite simple nonabelian group. 
Then $G=(A_5)^\N$ has no subgroups of countably infinite index.

\msk

This follows from Theorem 1.10 of Thomas [\thomas]. In fact Thomas
classifies infinite products of finite groups in which every subgroup
of index $<2^{\aleph_0}$ is necessarily open; such groups do not have
countably infinite index subgroups.
From Corollary 1.3 it follows that $G$ in Example 3.1 
has nonmeasurable subgroups. 
\msk

The literature provides some guidance on the situation of 
finite index subgroups.

\ssk

We begin with a result of
{\smc  M. G. Smith}  and {\smc J. S. Wilson} [\smith].
\msk\bf
Proposition 3.2. \quad \it Let $G$ be a profinite group. Then
all finite index normal subgroups are open if and only if
there are countably many finite index subgroups.
\qed

Since the cardinality of the set open normal subgroups in a profinite
group is the weight of the group, an immediate corollary is

\msk\bf
Corollary 3.3. \quad\it Let $G$ be a profinite group of uncountable weight.
Then $G$ contains  nonclosed  finite index subgroups  and these
are, accordingly, nonmeasurable.\qed  \rm

\msk
In fact as a consequence of {\smc Peterson}'s Theorem 1.2(2) of [\pete] it
has been known for some time that these large profinite groups contain
at least $2^{2^{\aleph_0}}$ nonmeasurable subgroups.

\msk

The focus therefore is on profinite groups with countably many
finite index normal subgroups, and in accordance with some authors
we use the following definition (see [\ribes], Section 4.2, pp. 124ff.)

\msk

\bf 
Definition 3.4. \quad \rm A {\it strongly complete group} is a 
profinite group 
in which every finite index  subgroup is open.
\msk

We summarize our findings:

\msk

\bf
Theorem 3.5.\quad \it A compact group in which every subgroup is
measurable is a strongly complete group. \bf

\msk

Proof. \quad \rm Let $G$ be a compact group. If it is not 
totally disconnected, then it has a subgroup of countably infinite
index by Theorem 2.3 and thus a nonmeasurable subgroup. 
 If all subgroups of $G$ are measurable,
 then $G$ is  profinite. If it has a 
subgroup of finite index that fails to be open closed, then
such a subgroup is nonmeasurable by Proposition 1.1(c). Thus
all finite index subgroups of $G$ are open closed and so $G$ 
is a strongly complete group. \qed

\msk
So Question 1 reduces to
 
\msk

{\bf Question 2.}  \quad {\it Does every infinite strongly complete 
 group contain a nonmeasurable subgroup?}

\msk


We keep in mind that  Smith and Wilson [\smith] showed that
a profinite group is strongly complete if and only if 
it has only countably
many finite index  subgroups. Such a group is necessarily metric.
Segal and Nikolov   [\segal] showed that all topologically finitely 
generated metric profinite groups are strongly complete,
as had been conjectured by Serre.
\msk 
Typical examples in this class of groups are countable 
products of pairwise nonisomorphic  simple finite groups.
A result of {\smc Saxl}'s and {\smc Wilson}'s [\saxl] says:

\msk\bf

Proposition 3.6.\quad\it  Let $\{G_n:n\in\N\}$ be a sequence of finite simple
nonabelian groups and $G=\prod_{n\in\N}G_n$.
Then the following conditions are equivalent:
{\parindent2em

\item{\rm(i)} Infinitely many of the $G_n$ are isomorphic.

\item{\rm(ii)} $G$ is not strongly complete.\qed 

}\rm

\noindent We have observed in Theorem 3.5 that (ii) implies

(iii) {\it G has nonmeasurable subgroups.}
 
\msk

\noindent The implication ``(i) implies (ii)'' also 
follows from our Corollary 1.4 above.

\bsk 

\cen{\bf 4. Metric compact groups}
\msk
\def\b #1{{\cal B}(#1)}    
\def\s #1{{\cal S}(#1)}    

We propose to calculate the cardinality of the set $\s G$
of (not necessarily closed!) subgroups of $G$.
We let $\cc=2^{\aleph_0}$ denote the cardinality of the
continuum and $\b G$ the set of all Borel subsets of $G$.
\msk\nin
\bf Proposition 4.1. \rm If $X$ is an infinite 2nd countable metric space,
then $$\card(\b X)\le\cc.$$

\ssk
\bf Proof.\quad \rm In [\bourb], Exercise 4 c),  \S6, Chap.~9
it is established
that the cardinality of the set of Borel subsets of a
metric second countable space is $\le\cc$. \qed

\nin
\bf Remark 4.2. \rm For Haar measure $\mu$ on a compact metric group $G$,
a subset $X$ is measurable iff there are sets $B_1, B_2\in \b G$ such that
$B_1\subseteq X\subseteq B_2$ such that $\mu(B_2\setminus X)=0=
\mu(X\setminus B_1)$.  \qed

\nin
(See e.g.\ [\rudinb],  10.10; the argument
given there is quite general.)

Since the cardinality of an infinite compact metric group is $\cc$ and
therefore the cardinality of the set of its subsets is $2^\cc$, the
following observation is trivial.

\msk\nin
\bf Lemma 4.3. \it If $G$ is a compact metric group, $G/N$  an infinite
quotient group, and $H$ any subgroup of $G$
 then
$\max\{\card(\s H), \card(\s {G/N})\}\le \card(\s G)\le 2^\cc$.\rm \qed

The following is easily established.
\msk\nin
\bf Lemma 4.4. \it A vector space of the infinite dimension $\aleph$
has $2^\aleph$ vector subspaces.

\bsk

We now prove that an infinite compact metric group has as
many subgroups as it has subsets.

\msk
\bf Theorem 4.5. \quad\it Let $G$ be an infinite metric  compact group. Then
$$\card(\s G)=2^\cc.$$
\msk
\bf Proof. \quad \rm  (i) By Zelmanov's Theorem [\zel], $G$ contains an
infinite abelian subgroup $A$ which we may assume to be closed. Then
$A$ is a compact metric abelian group. If the assertion of the Theorem
is true for abelian groups, then it is true in general by Lemma 4.3
Thus we assume from here on that $G$ is abelian.

\ssk

(ii) By Lemma 4.3, we have 
$$\max\{\card(\s {G_0}),\card(\s {G/G_0})\}\le
\card(\s G)\}.$$  If $G_0\ne\{0\}$, by Lemma 4.3, we may assume $G$ to be
connected. Then $G$ is divisible and the torsion free summand is
a $\Q$-vector space of dimension $\cc$. Hence $G$ has at least $2^\cc$
subgroups. So from here on we may and will assume that $G$ is profinite,

\ssk

(iii) For any prime $p$ Let $G_p$ be the $p$ primary component of
the abelian profinite group $G$.
Then $G\cong \prod_{p\ {\rm prime}}G_p$
[\compbook], Corollary 8.8(iii).
There are two cases:

\ssk

(A) There is a $p$ such that $G_p$ is infinite, or

\ssk

(B) for all $p$ the group $G_p$ is finite.

\ssk

(iv) In Case (A) we invoke Lemma 3.4 and assume that $G$ is a pro-$p$-group.
The $p$-socle of $\hat G$ is a GF$(p)$-vector space of dimension $p$-rank $G$
(see [\compbook] A1,21).

Then again we have two cases: (a) The $p$-rank, that is the
GF$(p)$-dimension of
the $p$-socle of $\widehat G$ is finite, or else, (b) it is infinite.

In case (a) $\widehat G$ has a direct summand isomorphic to $\Z_p$
and by Lemma 3.4 again, we may assume that $G=\Z_p$.
In this case $G$ has a subgroup which is isomorphic to a
$\Q$-vector space of dimenstion $\cc$ and thus by Lemma 4.4, $G$
has at least $2^\cc$ subgroups.

In case (b), by duality, there is
a surjective morphism $G\to \Z(p)^\N$. Then by Lemma 3.4
 we assume that
$G=\Z(p)^\N$. Then $G$ is algebraically a GF$(p)$-vector space of 
dimension $\cc$.
By Lemma 4.4 again, $G$ has at least $2^\cc$ subgroups.
\ssk

(v) We are now reduced to case (B) and so, since $G$ is infinite, there are
infinitely many primes $p_1, p_2,\dots$ such that $G_{p_n}\ne\{0\}$,
$G_{p_n}$ finite and we may assume w.l.o.g. $G_{p_n}=\Z(p_n)$.
So we have to prove that for $G=\prod_{n=1}^\infty \Z(p_n)$
we have $\card(\s G)=2^\cc)$.
Now $\bigoplus_{n=1}^\infty Z(p_n)$ is the torsion group $T$ of $G$,
and  we claim that the torsion free rank of $G/T$ is $2^\cc$.
This claim is equivalent to saying that $\card G/T=2^\cc$.
But since $\card G=2^\cc$ and $T$ is countable, this is clear. \qed
\msk

\bf Corollary 4.6. \quad  \it Every infinite
compact group has a subgroup which is
not a Borel subset.

\msk
\bf Proof.\quad \rm By Proposition 4.1 and Theorem 4.5, 
every infinite compact {\it metric group} has more 
subgroups than it has Borel subgroups.

By the results preceding Section 4, every compact group which fails to
be a metric profinite group has a subgroup which is nonmeasurable
for Haar measure. Since all Borel subgroups are Haar measurable,
none of these is a Borel subgroup. \qed

\bsk
\rm

\cen{\bf References}

\msk

\def\flqq{\raise1pt\hbox{$\scriptscriptstyle\ll$}}
\def\frqq{\raise.8pt\hbox{$\scriptscriptstyle\gg$}}

{\parindent 1.5em
\item{[\bourb]} Bourbaki, N., \flqq Topologie g\'en\'erale\frqq, Chap 9,
Hermann, Paris, 1958.

\item{[\bour]} Bourbaki, N.,  \flqq Int\'egration\frqq, Chap. 1,2,3 et 4,
Hermann, Paris, 2me \'ed, 1965.

\item{[\comf]} Comfort, W. W.,  S. U. Raczkowski, and F. J. Trigos-Arrieta,
{\it Making group topologies with, and without, convergent sequences},
Appl. Gen. Topol. {\bf7} (2006), 109--124.

\item{[\ferher]} Ferrer, Maria V., Hern\'andez, Salvador,
{\it Dual topologies on
non-abelian groups}, 
Topology and its Applications, Special Issue in Honor of Dikran Dikranjan
(2012), to appear.

\item{[\gart]}
Gartside, P., and M. Smith,
{\it Counting the closed subgroups of profinite groups},
{\tt arXiv:0809.4734v1 [math.GR] 27 Sep 2008}

\item{[\gelb]} Gelbaum, B. R.,
``Problems in real and complex analysis,''  Problem Books in Analysis,
Springer-Verlag, New York etc., 1992, x+488p. 

\item{[\cont]} Gierz, G., K. H. Hofmann, K. Keimel, J. D. Lawson,
M. W. Mislove, and D. S. Scott, ``Continuous Lattices and Domains,''
Encycl. of Math. and its Appl. {\bf93}, Cambridge Univ. Press, 2003,
xxxvi+591 pp.

\item{[\hart]}
Hart, J. E., and K. Kunen,
{\it Compactifications of Non-Abelian Groups},
Topology Proceedings 26 (2001-2002), 593--626

\item{[\hern]}
Hern\'andez, S., {\it Questions raised at the
conference} \rm ``Algebra meets Topology,'' Barcelona, July 2010.
{\tt hernande@mat.uji.es}

\item{[\HHM]} Hern\'andez, S., K. H. Hofmann, and S. A. Morris,
{\it The weights of closed subgroups of a locally compact group},
J. of Group Theory {\bf15} (2014), 613--630. 

\item{[\hewross]}
Hewitt, E., and K. A. Ross,
``Abstract Harmonic Analysis I,''
Grundlehren {\bf115} Springer Verlag, Berlin etc., 1963.

\item{[\hof]} Hofmann, K. H.,
Arc components of locally compact groups are
Borel sets, Bull. Austr. Math. Soc. {\bf65} (2002), 1--8

\item{[\layer]}  Hofmann, K. H., and S. A. Morris,
{\it A structure theorem on compact groups}, Math. Proc.
Camb. Phil. Soc. {\bf130} (2001), 409--426.

\item{[\compbook]} ---,
``The Structure of Compact Groups,'' de Gruyter Studies 
in Math. {\bf25}, Berlin, 3$^{\rm rd}$
Edition 2013,  xxii+924 pp.

\item{[\kall]} Kallman, R. R.,
{\it Every reasonably sized matrix group is a subgroup of $S\infty$,}
Fund. Math. 164 (2000), 35--40.

\item{[\khara]} Kharazishvili, A. B.,
{\it On thick subgroups of uncountable $\sigma$-compact locally compact
commutative groups},
Topology and its Applications {\bf156} (2009), 2364--2369.

\item{[\klepp]} Kleppner, A.,
{\it Measurable homomorphisms of locally compact groups,}
Proc. Amer. Math. Soc. {\bf106} (1989), 391--395.

\item{[\mont]} Montgomery, D. and L. Zippin,
``Topological Transformation Groups,''
Interscience, New York, 1955.

\item{[\segal]} Nikolov, N., and D. Segal,
{\it On finitely, generated profinite groups, 
I: Strong completeness and uniform bounds}
Ann. of Math. {\bf165} (2007), 171--238.

\item{[\pete]} Peterson, H. L.,
{\it Extensions of Haar measure to relatively large
nonmeasurable subgroups},
Trans. Amer. Math. Soc. {\bf228} (1977), 359--370.

\item{[\ribes]} Ribes, L., and P. Zaleskii,
``Profinite Groups,''
Springer-Verlag, Berlin etc., 2000, xiv+435 pp.
2nd edition 2010, xvi+464 pp.

\item{[\rudina]}
Rudin, W.,
``Fourier Analysis on Groups,''
Interscience, New York, 1962.

\item{[\rudinb]}
---,
  ``Principles of Mathematical Anaysis,''
MacGraw- Hill, New York etc., 1953, 1964.

\item{[\saeki]}
S. Saeki and Karl Stromberg,
{\it Measurable subgroups and nonmeasurable characters},
Math.Scnad. {\bf 57} (1985), 359--374.

\item{[\saxl]}
Saxl J., and J. S. Wilson,
{\it A note on popers in simple groups}
Math. Proc. Cambridge Phil. Soc. {\bf122} (1997), 91--94

\item{[\smith]}
Smith, M. G., and J. S. Wilson
{\it On subgroups of finite index in compact Hausdorff groups},
 Archiv d. Math. {\bf80} (2003), 123--129.

\item{[\thomas]} Thomas, S.,
{\it Infinite products of finite simple groups II},
J. of Group Theory {\bf2} (1999) 401--434.

\item{[\zele]} Zelenyuk, Yevhen,
``Ultrafilters and Topologies on Groups,''
DeGruyter Expositions in Math {\bf50}, Berlin, 2011

\item{[\zel]} Zelmanov, E. I.,
{\it On periodic compact groups}, Israel J. of Math. {\bf77} (1992), 83--95.

}

\bsk\noindent
{\smc Authors' Adresses}
\msk
\parindent=0pt

{Salvador Hern\'andez}

{hernande@mat.uji.es}

{Universitat Jaume I, INIT and Depto de  Matem\'aticas,
Campus de Riu Sec, 12071 Castell\'on, Spain}
\msk
{Karl H. Hofmann, Corresponding author}

{hofmann@mathematik.tu-darmstadt.de}

{Fachbereich Mathematik, Technische Universit\"at Darmstadt, Schlossgartenstrasse 7, 64289 Darmstadt, Germany. }

 \msk
{Sidney A. Morris}

{morris.sidney@gmail.com}

{School of Science, IT, and Engineering, Federation University Australia,
Victoria 3353, Australia, and School of Engineering and Mathematical Sciences,
La Trobe University, Bundoora, Victoria 3086, Australia}

\bye